\documentclass[12pt]{amsart}
\usepackage{amsfonts,amssymb,amsxtra}
\usepackage{latexsym}

\newcommand{\el}{\par \mbox{} \par \vspace{-0.5\baselineskip}}

\newcounter{amoi}
\newtheorem{theo}{Theorem}
\newtheorem{prop}{Proposition}
\newtheorem{coro}{Corollary}
\newtheorem{defi}{Definition}
\newtheorem{rk}{Remark}
\newtheorem{nota}{Notation}
\newtheorem{example}{Example}

\newcommand{\noi}{\noindent}

\newenvironment{pf}{\noi {\el \noi \bf Proof.}}{\hfill $\Box$ \el}

\newcommand{\Ocal}{\mathcal O}

  \usepackage{showlabels}

\begin{document}

\begin{center}{\Large{Modernity of Halphen's variations on a theme of Monge}\vskip 1cm} Caroline {\sc Gruson} {\footnote{Equipe de g\'eom\'etrie, U.M.R. 7502 du CNRS, Institut Elie Cartan, 
Universit\'e de Lorraine, BP 239,
54506 Vandoeuvre-les-Nancy Cedex, France. E-mail:
Caroline.Gruson@univ-lorraine.fr}}, Laurent {\sc Gruson} {\footnote
{L.M.V., Universit\'e de Versailles,  E-mail: Laurent.Gruson@uvsq.fr}}.

 \end{center}
\vskip 2cm
\begin{center}
{\it Dedicated to Philippe Ellia.}
\end{center}

\section{Introduction}

Among the mathematical researches of Georges-Henri Halphen (1844-1889) there is an almost unexplored {\it corpus} on the theme of infinitely near points in a projective space of dimension $2$ and $3$.
\el

{\noi}{\bf History.}

Halphen was an artillery officer educated at Ecole polytechnique (Paris). His professional life as a mathematician took place in 1872-1886: he was then on secondment as {\it r\'ep\'etiteur d'analyse}, Ecole polytechnique. This period was preceded by a formidable {\it amateur} work on problems in algebraic geometry, notably the beginnings of enumerative geometry and the classification of algebraic space curves; it was followed, after his reintegration in the army, by the elaboration of his unfinished {\it Trait\'e des fonctions elliptiques et de leurs applications}.
\el

It is easy to describe in modern terms the sequence of his main papers on infinitely near points:
\el

{\noi}he anticipated in \cite{H1} the idea of Nash blowing-up, as a process of desingularisation of plane curves (it may have been the first occurrence of this idea),
\el

{\noi}he initiated in \cite{H2} an enumerative theory of the variety of infinitely near points in the plane, of an arbitrary given order. This theory (which was probably conceived at the same time as his (now unearthed) theory of conics) is parallel to the theory of {\it Semple towers} \cite{Se} - it would be interesting to decide whether they are identical or not,
\el

{\noi}he initiated in \cite{H3} an invariant theory for the action of the projective group on the same variety. He extended this theory to the projective space in \cite{H4},
\el

{\noi}on the occasion of the competition for the {\it Grand prix des sciences math\'ematiques} (Acad\'emie des sciences, 1880) he extended this invariant theory to the $n$-dimensional projective space and identified it to the invariant theory of linear differential equations of order $n+1$ \cite{H5},
\el

{\noi}he began in \cite{H6} an incursion in the Galois theory of differential equations, on the occasion of a result by Goursat on a case of reduction of the fourth-order linear differential equation towards two second-order differential equations: he was able to describe an invariant of the given {4}-th order equation whose vanishing implies that the reduction is possible.
\el

{\noi}{\bf A long-standing ambition.}

The papers \cite{H3}, \cite{H4}, \cite{H5} are very closely related (the variation of context hides this at first glance). Let us try to extract their common governing idea in the language of the more developed \cite{H5}:
\el

{\noi}{\bf Problem.} An ordinary linear differential equation of order $n+1$ (in the complex domain) is susceptible to be transformed by 1) multiplying the unknown function by a given function, 2) substituting to the variable a given function of another variable. Is it possible to transform it into another equation belonging to one of the established "integrable" types (at the time of \cite{H5})? 
\el

This datum provides a group (the allowed transformations) acting on a vector space (the set of coefficients of the equation), the problem is to recognise the elements of some specified orbits. Halphen's approach is to describe the invariant functions of the equation. The first obvious task is to identify the {\it unstable} equations on which all invariants vanish identically: they are the Veronese transforms (the elevation to the power $n$ of the unknown function) of second-order equations, leading to the uncomfortable conclusion that any attempt of this type will be ineffective on second-order equations.

Halphen uses {\it inter alia} the idea of a {\it slice} (or {\it canonical form}): let us recall
\begin{defi} Let $G$ be a group acting on an algebraic variety $X$, $Y$ be a locally closed subvariety of $X$, $H$ be a subgroup of $G$ stabilising $Y$. The couple $(Y,H)$ is a slice for the action of $G$ on $X$ if a general orbit of $X$ under $G$ meets $Y$ along an orbit under $H$. \end{defi}

The construction of this slice is the goal of the present paper; here we are content with a short description. A linear differential equation of order $n+1$ gives rise, {\it via} its vector space of solutions, to an analytic curve $\mathcal{C}$ in the projective space $\mathbb{P}_{n}$. Halphen defines a section of $\omega_{\mathcal{C}}^{\otimes 3}$ which, unfortunately, should be denoted $\frac{\mu}{U^{3}}$ (where $U$ is the wronskian): 
\el

{\noi}for any point $O$ of $\mathcal{C}$ there is a unique normal rational curve $\Gamma_{O}$ of degree $n$ containing the thickening of $O$ on $\mathcal{C}$ of length $n+3$; the meaning of the relation $\mu=0$ is that the Zariski tangent plane to the double point $O$ of the curve $\mathcal{C}\cup\Gamma_{O}$ is contained in the osculating hyperplane to $(\mathcal{C},O)$.
\el 

The function $\mu$ of the germ $(\mathcal{C},O)$ (parametrised by the variable $t$ of the differential equation) is called (by Sylvester \cite{Sy}) the {\it Monge invariant}.

Assuming $\mu\neq 0$, Halphen defines an intrinsic {\it curvilinear abscissa} $s$ by $\frac{ds}{dt}=\frac{\mu^{\frac{1}{3}}}{U}$ and then, taking $s$ as the new variable, an intrinsic transform of the given differential equation. The coefficients of the transform, apart from the adjunction of a cubic root of $\mu$, are rational functions of the original coefficients.

Let us get back to the problem. The orbits considered are represented by differential equations whose coefficients are rational functions over a rational or elliptic base curve; the set of singular points of such an equation belongs to a specified type: for instance, following Schwarz, if the base curve is rational, it consists of at most three regular singular points whose exponents $p,q,r$ are integers satisfying $\frac{1}{p}+\frac{1}{q}+\frac{1}{r}\geq 1$. If we know the intrinsic transform of our equation, we can decide if the algebraic relations between the coefficients of this transform allow to rationally define another transform over a function field of transcendence degree $1$ and genus $0$ or $1$, and look at its singular points. 

A typical example is given in \cite{H7}: a family of third-order linear differential equations which may be defined over $\mathbb{P}_{1}$ with three singular points of unorthodox exponents $\frac{1}{2},\frac{1}{3},\frac{1}{7}$. 
\el

{\noi}{\bf Acknowledgements.} The second author is grateful to Abdelghani El Mazouni, who wrote \cite{EM} (an analysis of Halphen's thesis \cite{H3}), for very useful discussions relevant to the subject of the present paper. 

\section{Preliminaries}

Let $\mathbb{P}$ be an $n$-dimensional complex projective space ($n\geq 2$). We denote by $A$ the vector space $H^{0}(\Ocal_{\mathbb{P}}(1))$.

We begin with a (well-known) formal definition, translating the idea of a smooth algebroid curve $\mathcal{C}$ in $\mathbb{P}$ centered at a point $O$ of $\mathbb{P}$:

\begin{defi} A {\it smooth arc} in $\mathbb{P}$ is a formal series

$$a=\sum_{i\geq 0}a_{i}t^{i}\in A\spcheck\otimes\mathbb{C}[[t]]$$

{\noi}such that $a_{0}\wedge a_{1}\neq 0$, considered modulo the (right) action of the group $G=G_{0}\rtimes G_{1}$, where $G_{0}$ is the multiplicative group of power series beginning with a non-zero constant term (acting by multiplication on $A\spcheck\otimes\mathbb{C}[[t]]$) and $G_{1}$ is the group consisting of power series beginning with a non-zero term of degree $1$ (under the composition law $(u,v)\mapsto v\circ u$), acting by composition on both $G_{0}$ and  $A\spcheck\otimes\mathbb{C}[[t]]$. \end{defi}

In what follows we freely use schemes over $\mathbb{C}$ which are not necessarily of finite type, the basic example being the vector space $\mathbb{C}[[t]]$, which we view as the spectrum of the symmetric algebra of its topological dual $\mathbb{C}((t))/\mathbb{C}[[t]]$.

\begin{rk} {\normalfont1. If we replace $\mathbb{C}[[t]]$ by the ring $\mathbb{C}\lbrace{t}\rbrace$ of convergent power series, this definition amounts to describing the set of germs of analytic curves $(\mathcal{C},O)$ in $\mathbb{P}$, ascribing to a small value of $t$ the point $[a(t)]\in\mathbb{P}$ and marking the point $O=[a(0)]$.

2. The elements of the (Zariski)-open set $\mathcal{P}\subset A\spcheck\otimes \mathbb{C}[[t]]$ defined by $a_{0}\wedge a_{1}\neq 0$ are called {\it parametrised smooth arcs}. Since the complement of $\mathcal{P}$ in $A\spcheck\otimes\mathbb{C}[[t]]$ is of codimension $n\geq 2$, the regular functions over $\mathcal{P}$ extend to $A\spcheck\otimes\mathbb{C}[[t]]$.

3. A truncated smooth arc of length $r+1$, or order $r$, is defined in a similar way: replace the algebra $\mathbb{C}[[t]]$ by its quotient $\mathbb{C}[[t]]/(t^{r+1})$.}
\end{rk} 

We denote by $\mathcal{A}_{\mathbb{P}}$, resp. $\mathcal{A}_{\mathbb{P},r}$, the set of smooth arcs, resp. truncated smooth arcs of order $r$, in $\mathbb{P}$; we will omit the subscript $\mathbb{P}$ when the context is clear.
\el

Let us construct local cross-sections for the action of $G$ on $\mathcal{P}$. Choose coordinates $X_{0},\ldots,X_{n}$ on $\mathbb{P}$ and fix two distinct integers $i,j$ in $[0,n]$. Denote $\mathcal{U}_{ij}$ the subset of $\mathcal{A}_{\mathbb{P}}$ consisting of smooth arcs $(\mathcal{C},O)$ such that the rational function $\frac{X_{j}}{X_{i}}$ on $\mathbb{P}$ is defined at $O$ and reversible on $\mathcal{C}$ at the point $O$. 

Take $\mathcal{U}_{01}$ for instance: we may identify it with the closed subset of $\mathbb{C}^{n+1}[[t]]$ of which a typical element is the $(n+1)$-uple $(1,x_{1}+t,\Sigma_{i\geq 0}b_{2,i}t^{i},\ldots,\Sigma_{i\geq 0}b_{n,i}t^{i})$.

Each $\mathcal{U}_{ij}$ is thus an infinite-dimensional affine space, and the change of coordinates from $\mathcal{U}_{ij}$ to $\mathcal{U}_{kl}$ is a rational map defined on their open intersection $\mathcal{U}_{ij,kl}$ (in $\mathcal{A}$!), with closed graph in $\mathcal{U}_{ij}\times\mathcal{U}_{kl}$. Therefore $\mathcal{A}$ becomes an infinite-dimensional variety, the projective limit on $r$ of the finite-dimensional varieties $\mathcal{A}_{r}$. The action of $G$ on $\mathcal{P}$ is {\it free}: the projection of its graph in $G\times\mathcal{P}\times\mathcal{P}$ into the factor $\mathcal{P}\times\mathcal{P}$ is a closed immersion; moreover $\mathcal{A}$ is the corresponding geometric quotient, cf.\cite{GIT}, def. 0.8 and 0.6. 
\el

{\noi}{\bf The Picard group of $\mathcal{A}$.}

The character group $X(G)$ of $G$ is the direct product $X(G_{0})\times X(G_{1})$. Both factors become isomorphic to $\mathbb{Z}$ by choosing the reduction character $$a\colon=\Sigma_{i\geq 0}a_{i}t^{i}\mapsto \bar{a}\colon=a_{0}$$ of $G_{0}$, resp. $$b\colon=\Sigma_{i\geq 1}b_{i}t^{i}\mapsto \bar{b}\colon=b_{1}$$ of $G_{1}$. We already noticed that the complement of $\mathcal{P}$ in the affine space $Spec(A\spcheck\otimes\mathbb{C}[[t]])$ is of codimension $n$; therefore, since $G$ acts freely on $\mathcal{P}$, the Picard group of $\mathcal{A}=\mathcal{P}/G$ is naturally isomorphic to the character group of $G$ (cf.\cite{GIT}, p.32): we associate to an element $(\delta,\varpi)\in\mathbb{Z}$ the line bundle $\Ocal_{\mathcal{A}}(\delta,\varpi)$ over $\mathcal{A}$, quotient of the trivial bundle $\mathcal{P}\times\mathbb{C}$ over $\mathcal{P}$ by the action of $G$ given on the factor $\mathbb{C}$ by multiplication by the character $(a,b)\mapsto\bar{a}^{\delta}.\bar{b}^{\varpi}$. The integer $\delta$, resp. $\varpi$, is called the {\it degree}, resp. {\it weight}, of the invertible sheaf $\Ocal_{\mathcal{A}}(\delta,\varpi)$. 

The vector space $H^{0}(\Ocal_{\mathcal{A}}(\delta,\varpi))$ is a subspace of $S{d}(A\otimes(\mathbb{C}((t))/\mathbb{C}[[t]]))$ (by construction); its elements will (seldom) be called {\it differential expressions} of degree $\delta$ and weight $\varpi$ over $\mathcal{A}$. We have $\mathcal{A}=\varprojlim\mathcal{A}_{r}$ and therefore $H^{0}(\Ocal_{\mathcal{A}}(\delta,\varpi))=\varinjlim H^{0}(\Ocal_{\mathcal{A}_{r}}(\delta,\varpi))$. Hence, given a differential expression $s\in H^{0}(\Ocal_{\mathcal{A}}(\delta,\varpi))$, there exists a smallest integer $r$ such that $s$ is defined over $\mathcal{A}_{r}$: it is called the $\it order$ of the expression $s$.

Taking again a system $X_{0},\ldots,X_{n}$ of homogeneous coordinates on $\mathbb{P}$, we may write an element $a(t)$ of $A\spcheck\otimes\mathbb{C}[[t]]$ as an $(n+1)\times\infty$-matrix $(a_{ij})_{i\in[[0,n]],j\in\mathbb{N}}$, considering that the $i$-th coordinate of $a$ is $\Sigma_{j\geq 0}a_{ij}t^{j}$. Then the function $a_{00}$, resp. $a_{00}a_{11}-a_{01}a_{10}$, is a section of $\Ocal_{\mathcal{A}}(1,0)$, resp. $\Ocal_{\mathcal{A}}(2,1)$. These sections do not vanish on $\mathcal{U}_{01}$. This allows to identify (by dividing by a suitable, non necessarily positive, monomial in these two particular sections) the differential expressions to unique polynomial expressions on $\mathcal{U}_{01}$. 

\begin{prop} Let $C\subset\mathbb{P}$ be a smooth curve. Let $\gamma:C\rightarrow\mathcal{A}$ be the morphism associating, to a point $p$ of $C$, the arc $(C,p)$. For $(\delta,\varpi)\in Pic(\mathcal{A})$ one has $\gamma^{*}(\Ocal(\delta,\varpi))=\omega^{\otimes\varpi}(\delta)$. \end{prop}
This is a boring exercise using local charts.

\begin{rk}\label{functions} {\normalfont Given a homogeneous polynomial function $f$ on $A\spcheck\otimes\mathbb{C}[[t]]$ (an element of $S^{\delta}(A\otimes(\mathbb{C}((t))/\mathbb{C}[[t]]))$) and a parametrised smooth arc $a\in A\spcheck\otimes\mathbb{C}[[t]]$, we form the Taylor series of $a$, $$a(t+h)=\sum_{i\geq 0}\frac{1}{i!}a^{(i)}(t)h^{i}$$ and evaluate $f$ on its coefficients. We thus obtain a regular function from $A\spcheck\otimes\mathbb{C}[[t]]$ to $\mathbb{C}[[t]]$, which we denote $f(t)$ and whose natural interpretation is the evaluation of $f$ at the same arc centered at the point image of $t$.

For instance, let $f$ be a section of $\Ocal_{\mathcal{A}}(\delta,\varpi)$ and $C\subset\mathbb{P}$ be a smooth curve (maybe analytic or even algebroid). For a given point $O$ of $\mathbb{C}$ we may represent the germ $(C,O)$ by a formal series $a(t)\in A\spcheck\otimes\mathbb{C}[[t]]$, this allows to define the generating elements (trivialisations) $a(t)$ of $\Ocal_{C}(-1)$ and $dt$ of $\omega_{C}$ (considered as invertible $\Ocal_{C}$-modules) at the point $O$. Then $f(t)$ is the evaluation at $O$ of the section of $\omega^{\otimes\varpi}(\delta)$ deduced from $f$ in this trivialisation.}\end{rk}

{\noi}{\bf Now enters the central object of this study, namely the differential invariants.}

\begin{defi}\label{defweight} A differential invariant of degree $\delta$ and weight $\varpi$ on $\mathbb{P}$ is a section of the invertible $\Ocal_{\mathcal{A}}$-module $\Ocal_{\mathcal{A}}(\delta,\varpi)$ which is invariant under the natural (left) action of the group $SL(A)$.\end{defi}

\begin{nota} {\normalfont The ring of differential invariants over $\mathbb{P}$ is denoted $\mathcal{H}_{n}$: note that this ring only depends on the dimension $n$ of the projective space. 

The degree of a differential invariant $I$, of bidegree $(\delta,\varpi)$, is always divisible by $n+1$; the quotient $d:=\frac{\delta}{n+1}$ is called the {\it reduced degree} of $I$. For reasons to be explained soon (cf. th.\ref{iota}), the integer $p\colon=\varpi-\binom{n+1}{2}d$ is called the {\it reduced weight} of $I$, and $(d,p)$ is called the {\it reduced bidegree} of $I$.}\end{nota}

\begin{example}\label{U}  The wronskian: {\normalfont the function $$a=\Sigma_{i\geq 0}a_{i}t^{i}\mapsto det(a_{0},\ldots,a_{n}): A\spcheck\otimes\mathbb{C}[[t]]\rightarrow\mathbb{C}$$ is a differential invariant of reduced bidegree $(1,0)$. We will follow Halphen in denoting it $U$.

This invariant vanishes exactly on inflexional arcs: more precisely, if $(\mathcal{C},O)$ is an arc represented by the series $a=\Sigma_{i\geq 0}a_{i}t^{i}$, define inductively the sequence $(i_{0},\ldots,i_{n})$ by setting $i_{j+1}$ as the smallest integer such that $a_{i_{j+1}}$ is not in the linear span of $(a_{i_{0}},\ldots,a_{i_{j}})$, or $+\infty$. The arc $(\mathcal{C},O)$ is called {\it inflexional} if and only the non-decreasing sequence $(i_{j}-j)_{j\in[0,n]}$ (sometimes called the Weierstrass gap sequence) is non-identically $0$: indeed, the integer $\Sigma_{j\in[0,n]}(i_{j}-j)$ is the order of the series $U(t)$.} \end{example}

\begin{example} \label{conics} {\normalfont First assume $n=2$, let $(\mathcal{C},O)$ be a smooth arc represented by the power series $a=\Sigma_{i\geq 0}a_{i}t^{i}\in A\spcheck\otimes\mathbb{C}[[t]]$. Consider its image by the second Veronese map $$v_{2}:\mathbb{P}\rightarrow P(S^{2}A\spcheck),$$ it is represented by the series $$\Sigma_{k\geq 0}(\Sigma_{i+j=k}a_{i}a_{j})t^{k}\in S^{2}A\spcheck\otimes\mathbb{C}[[t]].$$ The wronskian of this arc is a new differential invariant which happens to be divisible by $U$. The quotient will be denoted $\mu$, its reduced bidegree is $(3,3)$.

Assume now $n=3$, let $(\mathcal{C},O)$ be a smooth arc represented by the power series $a=\Sigma_{i\geq 0}a_{i}t^{i}\in A\spcheck\otimes\mathbb{C}[[t]]$. The {\it tangent scroll} $(T_{\mathcal{C}},t_{O})$ of $(\mathcal{C},O)$ (i.e. the locus of tangents to $\mathcal{C}$) may be represented by the power series $$\Sigma_{k\geq 0}(\Sigma_{i+j=k}(i+1)a_{i+1}\wedge a_{j})t^{k}\in\wedge^{2}A\spcheck\otimes\mathbb{C}[[t]].$$ The wronskian of this arc is a new differential invariant of reduced bidegree $(3,3)$. We will denote it $\mu$.

These two invariants will be unified later as the {\it Monge invariant}, defined in any dimension $n$.}\end{example}

Let $\mathcal{G}_{n+1}$ denote the grassmannian $Grass(n+1,\mathbb{C}[[t]])=\varprojlim_{r} Grass(n+1,\mathbb{C}[t]/(t^{r}))$ (viewed as a pro-algebraic variety). The torus $\mathbb{C}^{*}$ acts on the ring $\mathbb{C}[[t]]$ (by $(a,t)\mapsto at$) and on $\mathcal{G}_{n+1}$, this action gives the coordinate ring of $\mathcal{G}_{n+1}$ a second graduation, called the weight in consistency with Definition \ref{defweight}.

\begin{nota}\label{partitions}{\normalfont For every partition $\lambda=(\lambda_{0}\geq\lambda_{1}\geq\ldots\geq\lambda_{n})$ of length $\leq n+1$, set $\vert\lambda\vert\colon=\lambda_{0}+\ldots+\lambda_{n}$ and

\begin{center} $U_{\lambda}(\Sigma_{i\geq 0}a_{i}t^{i})\colon=\det(a_{\lambda_{n}},a_{\lambda_{n-1}+1},\ldots,a_{\lambda_{0}+n})$:\end{center}

{\noi}this is a well-known section of $\mathcal{O}_{\mathcal{G}_{n+1}}(1)$ (Pl\"{u}cker coordinate with respect to the topological basis $(t^{i})_{i\in\mathbb{N}}$), of bidegree $(1,\binom{n+1}{2}+\vert\lambda\vert)$: by a slight abuse we will say that its reduced bidegree is $(1,\vert\lambda\vert)$.}
\end{nota}

\begin{prop} Every differential invariant on $\mathbb{P}$, of reduced bidegree $(d,p)$, can be written as a bihomogeneous polynomial of degree $(d,p)$ in the $U_{\lambda}$, assuming the preceding convention for the reduced degree of $U_{\lambda}$. This polynomial is uniquely determined modulo the Pl\"{u}cker relations between the $U_{\lambda}$.\end{prop}

This is clear once it is noted that the quotient of $A\spcheck\otimes\mathbb{C}[[t]]$ under $GL(A\spcheck)$ is the grassmannian $\mathcal{G}_{n+1}\colon=Grass(n+1,\mathbb{C}[[t]])$.

\begin{rk}\normalfont {The set of non-inflexional smooth parametrised arcs in $\mathbb{P}$ (considered modulo $SL(A\spcheck)$) is the affine open set of $\mathcal{G}_{n+1}$ defined by $U\neq 0$; the Pl\"{u}cker coordinates $u_{i,k}\colon=\frac{U_{\lambda_{i,k}}}{U_{0}}$, where $k\geq 1$ and $\lambda_{i,k}$ is the partition $(k,1,\ldots,1,0,\ldots,0)$ of length $i\in[1,n+1]$ of $k+i-1$, are the natural coordinates of this affine variety, whose coordinate ring $\mathbb{C}[u_{ik}]$ remains $\mathbb{Z}$-graded by the reduced weight: $\deg(u_{ik})=k+i-1$. The differential invariants on $\mathbb{P}$ can (after multiplication by a suitable power of $U$) be written in a unique way as homogeneous polynomial functions of these variables.

In his paper \cite{H4} about the case $n=3$, Halphen systematically writes his invariants according to the convention $U=1$.} \end{rk}    

{\noi}{\bf An inductive construction of differential invariants.}

\begin{theo}\label{iota} There is an injective homomorphism of $\mathbb{C}$-algebras $\iota_{n}:\mathcal{H}_{n}\rightarrow\mathcal{H}_{n+1}$, characterised by the following property:

for any $I\in\mathcal{H}_{n}$, $\iota_{n}(I)$ vanishes on the germ $(\mathcal{C},O)$ if and only if $I$ vanishes on the projected germ of $(\mathcal{C},O)$ from the point $O$. 

Moreover $\iota_{n}$ preserves the reduced bidegree of differential invariants.\end{theo} 
\begin{pf} For a smooth arc $(\mathcal{C},O)$ in $\mathbb{P}_{n+1}$ represented by 
$a=\Sigma_{i\in\mathbb{N}}a_{i}t^{i}\in A\spcheck\otimes\mathbb{C}[[t]]$, the projected arc from $O$ in the "hyperplane of infinity" (the set $\mathbb{P}_{O}$ of lines through $O$, or any hyperplane of $\mathbb{P}$ not containing $O$) may be represented by the series $\Sigma_{i\geq 1}\bar{a}_{i}t^{i-1}\in (A\spcheck/\mathbb{C}a_{0})\otimes\mathbb{C}[[t]]$ (where the class of $a_{i}$ modulo $\mathbb{C}a_{0}$ is denoted $\bar{a}_{i}$). It is smooth under the assumption \begin{equation}\label{gap}a_{0}\wedge a_{1}\wedge a_{2}\neq 0 \end{equation}
{\noi}(notice that the codimension of the closed subset of $\mathcal{A}$ defined by $a_{0}\wedge a_{1}\wedge a_{2}=0$ is $n-1$.)

Consider the local chart $\mathcal{U}_{01}$ of $\mathcal{A}_{\mathbb{P}_{n+1}}$. We write the coordinate functions on $\mathcal{U}_{01}$ as an $(n+2)\times\infty$-matrix $M$ (the entry $X_{ij}$ of $M$ is the coefficient of $t^{j}$ in the $i$-th coordinate series) which takes the form 
\begin{equation}\label{arcmatrix}
M=\left(\begin{array}{lcccccccccr}1&0&0&\ldots\\X_{1,0}&1&0&\ldots\\X_{2,0}&X_{2,1}&X_{2,2}&\ldots\\\vdots&\vdots&\vdots\\X_{n+1,0}&X_{n+1,1}&X_{n+1,2}&\ldots\end{array}\right).
\end{equation}
 Let $M_{00}$ be the $(n+1)\times\infty$-matrix obtained in deleting the $0$-th line and column of $M$. 

A differential invariant $I$ on $\mathbb{P}_{n}$ of reduced bidegree $(d,p)$ is in particular a section of $\Ocal_{\mathcal{G}_{n+1}}(d)$, therefore it can be evaluated at the matrix $M_{00}$: the result is a polynomial function on $\mathcal{U}_{01}$. The meaning of the vanishing of this function at $(\mathcal{C},O)$ is (under the condition (\ref{gap})) the vanishing of $I$ at the projection of $(\mathcal{C},O)$ from $O$. Therefore the closure in $\mathcal{A}_{\mathbb{P}_{n+1}}$ of the hypersurface $I(M_{00})=0$ of $\mathcal{U}_{01}$ is independent of any choice of coordinates, i.e. its equation is a differential invariant which we take as $\iota_{n}(I)$. Its evaluation over $\mathcal{U}_{01}$ shows that its reduced bidegree is $(d,p)$.  \end{pf}

{\noi}{\bf Link with differential equations.}
\el

Let $a\in A\spcheck\otimes\mathbb{C}[[t]]$ (with $\dim(A)=n+1$). Consider the sequence
\begin{center} $a, a^{[1]},\ldots,a^{[n+1]}$ \end{center}
where $a^{[i]}$ denotes the {\it divided derivative} $\frac{a^{(i)}}{i!}$. The Cramer formula implies the following relation
\begin{equation}\label{PV} \Sigma_{i\in[0,n+1]}(-1)^{i}a^{[i]}\det(a,\ldots,\hat{a}^{[i]},\ldots,a^{[n+1]})=0. \end{equation}
In the equation (\ref{PV}), the coefficient of $a^{[i]}$ is immediately identified with $U_{(1^{n+1-i})}(t)$ (the index is the partition consisting of $(n+1-i)$ parts all equal to $1$). Following the useful normalisation of coefficients made by Halphen, we get

\begin{prop}\label{equadif} Let $(\mathcal{C},O)$ be a non-inflexional arc represented by the power series $a$. Then $a$ satisfies the differential equation of order $n+1$
\begin{center} $a^{(n+1)}(t)+\Sigma_{1\leq i\leq n+1}\binom{n+1}{i}p_{i}(t)a^{(n+1-i)}(t)=0$\end{center}
{\noi}where $p_{i}(t)$ is the series $(-1)^{i}i!\frac{U_{(1^{i}) }(t)}{U(t)}\in\mathbb{C}[[t]]$. \end{prop}

{\noi}(Since we have assumed that $(\mathcal{C},O)$ is non-inflexional, we may invert the wronskian $U(t)$ in $\mathbb{C}[[t]]$.) Conversely, given a linear differential equation \begin{center}$x^{(n+1)}(t)+\Sigma_{1\leq i\leq n+1}\binom{n+1}{i}p_{i}(t)x^{(n+1-i)}(t)=0$\end{center} {\noi}of order $n+1$, with coefficients in $\mathbb{C}[[t]]$, its set of solutions $\mathcal{S}$ is an $(n+1)$-dimensional subvector space of $\mathbb{C}[[t]]$; the associated element of $\mathcal{S}\spcheck\otimes\mathbb{C}[[t]]$ is a non-inflexional arc. These two constructions are mutually reciprocal, a well-known fact since the 18-th century.
\el

Halphen uses rather insistently the correlation ({\it via} the link above) between duality of curves and adjunction of linear differential equations. Let us follow his first steps:

The Weyl algebra $\mathcal{W}$ of differential operators on $\mathbb{C}[[t]]$ is generated over $\mathbb{C}[[t]]$ by the element $\frac{d}{dt}$, subject to the relation $[\frac{d}{dt},t]=1$ (of course $[,]$ is the Lie bracket; an element $a$ of $\mathbb{C}[[t]]$ is viewed in $\mathcal{W}$ as the multiplication $m_{a}$ by $a$). The classical definition of the adjunction $L\mapsto L^{*}$ is the anti-involution of $\mathcal{W}$ fixing $m_{t}$ and transforming $\frac{d}{dt}$ into $-\frac{d}{dt}$. The adjoint of the linear differential operator $L$ defined by
\begin{center} $Lu\colon=u^{(n)}+p_{1}u^{(n-1)}+\ldots+p_{n}u$ \end{center}
{\noi}is the linear differential operator $L^{*}$ defined by
\begin{center} $L^{*}u\colon=(-1)^{n}u^{(n)}+(-1)^{(n-1)}\frac{d^{n-1}}{dt^{n-1}}(p_{1}u)+\ldots+p_{n}u$ \end{center}
{\noi}which may be written $(-1)^{n}(u^{(n)}+q_{1}u^{(n-1)}+\ldots+q_{n}u)$, with
\begin{center} $q_{i}=\Sigma_{0\leq k\leq n-i}(-1)^{i+k}\binom{n-i-k}{k}p_{i+k}^{(k)}.$\end{center}
Notice the important identity
\begin{center} $(Lu).v=u.L^{*}v+\frac{d}{dt}(B(u,v))$ \end{center}
{\noi}where $B(u,v)$ is the bilinear expression in $(u,v)$ given by
\begin{equation}\label{duality} B(u,v)\colon=\Sigma_{k,l}u^{(k)}v^{(l)}\Sigma_{j}(-1)^{j}\binom{j}{k}p_{n-l-j-1}^{(j-k)}.\end{equation}

The following proposition is classical (and essentially found in Bourbaki, {\it Fonctions d'une variable r\'eelle} IV, \S 2, prop. 6):
\begin{prop}\label{duality1} If $\mathcal{S}$, resp. $\mathcal{S}^{*}$, is the space of solutions in $\mathbb{C}[[t]]$ of the (above) linear differential operator $L$, resp. of its adjoint $L^{*}$, the bilinear form $(u,v)\mapsto B(u,v)(0)$ on $\mathcal{S}\times\mathcal{S}^{*}$ is non-degenerate. Identify each of these spaces to the dual of the other. The corresponding non-inflexional arcs in $\mathbb{P}(\mathcal{S}^{*})$ and $\mathbb{P}(\mathcal{S})$ are exchanged by the Pl\"{u}cker duality.\end{prop}

\section{The normal rational curve associated to a non-inflexional arc.}

Recall the classical fact that a group of $n+3$ distinct points in $\mathbb{P}$ in general linear position is contained in a unique normal rational curve of degree $n$. Eisenbud and Harris \cite{EH} extended this result to non-reduced finite subschemes of $\mathbb{P}$. We need a particular case of it:

Let $(\mathcal{C},O)$ be a non-inflexional arc in $\mathbb{P}$, defined by a power series $a\in A\spcheck\otimes\mathbb{C}[[t]]$. Let $R$ be the quotient ring $\mathbb{C}[[t]]/(t^{n+3})$. The linear map $A\otimes R\rightarrow R$ induced by $a$ defines a morphism $spec(R)\rightarrow\mathbb{P} $ which is a closed immersion (since the arc is smooth): this scheme is the {\it thick point} of order $n+2$, or length $n+3$, on the arc $(\mathcal{C},O)$.

\begin{theo}\label{Eisenbud-Harris}\cite{EH} The scheme $spec(R)$ is contained in a unique normal rational curve $\Gamma$ of degree $n$ of $\mathbb{P}$.\end{theo}

Since $\Gamma$ is rational, its Picard group has a unique ample generator: let $W$ be its two-dimensional space of sections. Then $\Gamma$ is the image of $\mathbb{P}(W\spcheck)$ by its $n$-th Veronese embedding which we denote $v_{n}$, and $A$ is thus identified with $S^{n}W$. We choose a basis $e_{0},e_{1}$ of $W$ such that $O=v_{n}([e\spcheck_{1}])$ (as usual the dual basis is denoted $e\spcheck_{0},e\spcheck_{1}$). We will use the basis $\epsilon_{i}=\binom{n}{i}e_{0}^{n-i}e_{1}^{i}$ of $A$, whose dual basis (up to the factor $n!$) is $\epsilon\spcheck_{i}=(e\spcheck_{0})^{n-i}(e\spcheck_{1})^{i}$.

Let $B$ be the stabiliser of $O$ in $SL(W\spcheck)$, a Borel subgroup. The following statement is a direct consequence of Theorem \ref{Eisenbud-Harris}:

\begin{coro} Denote $\mathcal{A}^{o}\subset\mathcal{A}$ the open subset of non-inflexional arcs. Let $\mathcal{S}_{B}$ be the set of arcs taking the form $(1,t,t^{2}+O(t^{n+3}),\ldots,t^{n}+O(t^{n+3}))$ in the coordinates $(\epsilon_{i})_{i\in[0,n]}$ of $\mathbb{P}$ (recall that $O(t^{i})$ means a power series of order $\geq i$). 

Then any $SL(A\spcheck)$-orbit in $\mathcal{A}^{o}$ meets $\mathcal{S}_{B}$ along a $B$-orbit, in particular $\mathcal{S}_{B}$ is a $B$-slice for the action of $SL(A\spcheck)$ on $\mathcal{A}$. 
\end{coro}

Indeed, the parametric representation of $\Gamma$ in the chosen coordinates is by construction $t\mapsto (1,t,\ldots,t^{n})$. If we use the quotient of the first two coordinates as the parameter $t$ for $(\mathcal{C},O)$ we get the corollary.

Let us be a little more precise about the action of the Borel subgroup $B$ of $SL(W\spcheck)$ in our situation: let $$b_{\alpha,\beta}=\left(\begin{array}{lccccr}\alpha&0\\\beta&\alpha^{-1}\end{array}\right)$$ be a typical element of $B$. Its action on $A\spcheck$ is given by the matrix 
$$M_{\alpha,\beta}=\left(\begin{array}{lccccccccccr}\alpha^{n}&0&0&\ldots&0\\\alpha^{n-1}\beta&\alpha^{n-2}&0&\ldots&0\\\alpha^{n-2}\beta^{2}&2\alpha^{n-3}\beta&\alpha^{n-4}&\ldots&0\\\vdots&\vdots&\vdots&&\vdots\\\beta^{n}&n\alpha^{-1}\beta^{n-1}&\binom{n}{2}\alpha^{-2}\beta^{n-2}&\ldots&\alpha^{-n}\end{array}\right).$$

We now define a new $SL(A\spcheck)$-invariant closed subset $\mathcal{M}$ of $\mathcal{A}^{o}$. Take a non-inflexional arc $(\mathcal{C},O)$ and the rational normal curve $\Gamma$ associated to it. The point $O$ of the reducible curve $\mathcal{C}\cup\Gamma$ clearly has multiplicity $2$, therefore the dimension of the Zariski tangent plane $P$ to this curve at $O$ is $2$ (Bourbaki, {\it Alg\`ebre commutative} VIII, \S7, exerc. 6). This plane, which we call the {\it principal plane} (associated to the pair $(\mathcal C, \Gamma)$),  contains the common tangent to $\mathcal{C}$ and $\Gamma$ and is contained in the tangent space to $\mathbb{P}$ at $O$; a hyperplane of $\mathbb{P}$, with equation $l$, contains $P$ if and only if $l(a(t)-(1,t,\ldots,t^{n}))=O(t^{n+4})$. 

We study the relative position of $P$ and the common osculating flag to $(\mathcal{C},O)$ and $\Gamma$: either 1) $P$ is not contained in the common osculating hyperplane to $\mathbb{C}$ and $\Gamma$ (the general case), or 2) it is: clearly this condition defines a closed, $SL(A\spcheck)$-invariant hypersurface $\mathcal{M}$ of $\mathcal{A}^{o}$. 

\begin{defi} We follow Sylvester \cite{Sy} (or at least we believe to do so) in calling $\mathcal{M}$ the {\it Monge hypersurface}. If needed, we will add a indexation by the ambiant projective space ($\mathcal M _\mathbb P$ here).\end{defi}

In case 1) the linear span of $P$ and the $(n-2)$-osculating linear variety to $(\Gamma,O)$ cuts $\Gamma$ in a point $\Omega\neq O$. The stabiliser $T$ of $\Omega$ in $B$ is a torus; if the basis of $W$ is modified so that $\Omega=[e\spcheck_{0}]$, the penultimate coordinate of $(\mathcal{C},O)$ becomes $t^{n-1}+O(t^{n+4})$, while the last one is $t^{n}+\alpha t^{n+3}+O(t^{n+4})$ for some non-zero scalar $\alpha$. Choose a sixth root $\alpha^{1/6}$ of $\alpha$ and the corresponding element $$\tau\colon=\left(\begin{array}{lccccr}\alpha^{1/6}&0\\0&\alpha^{-1/6}\end{array}\right)$$ {\noi}of $T$. After changing the basis according to $\tau$, the last coordinate is $t^{n}+t^{n+3}+O(t^{n+4})$. Therefore we get

\begin{theo}\label{3-slice} Let $\mu_{3}$ be the subgroup (of order $3$) of $PSL_{2}$ generated by the class of $$\left(\begin{array}{lccccr}e^{i\pi/3}&0\\0&e^{-i\pi/3}\end{array}\right).$$ Identify $A$ and $S^{n}(\mathbb{C}^{2})$ as above (cf. Theorem \ref{Eisenbud-Harris}) and set $\mathcal{S}_{\mu_{3}}$ as the closed subset of $\mathcal{U}_{01}$ whose typical element is 
\begin{center}$(1,t,t^{2}+O(t^{n+3}),\ldots,t^{n-2}+O(t^{n+3}),t^{n-1}+O(t^{n+4}),t^{n}+t^{n+3}+O(t^{n+4}))$\end{center} {\noi}if $n\geq 3$, resp. \begin{center}$(1,t,t^{2}+t^{5}+O(t^{7}))$\end{center} {\noi}if $n=2$. Then each $SL(A\spcheck)$-orbit in $\mathcal{A}^{o}\setminus\mathcal{M}$ meets $\mathcal{S}_{\mu_{3}}$ along a $\mu_{3}$-orbit. In particular $\mathcal{S}_{\mu_{3}}$ is a $\mu_{3}$-slice of $\mathcal{A}$. 
\end{theo}
The special case $n=2$ is examined in \cite{EM}, lemme 5.1.

\begin{rk}\label{n+4} The case 2) is characterised by the existence of a basis in which the coordinates of the arc are $\left(1,t, t^2+O(t^{n+3}), \ldots, t^{n-1}+O(t^{n+3}), t^n+O(t^{n+4})\right)$. Such a basis is not unique.
\end{rk}

\section{The Monge equation.}

Let us recall that the Monge hypersurface $\mathcal M$ of $\mathcal A$ is the closure of the set of non-inflexional arcs $(\mathcal C , O)$ such that the principal plane $P$ is contained in the osculating hyperplane.

\begin{theo}\label{Monge} The equation of $\mathcal M$ is
\begin{equation*}U^2 (U_3+2U_{21}-2U_{111})-3UU_1(U_{11}+U_{2})+2U_{1}^{3}=0,\end{equation*} 
where we use Notation \ref{partitions}.
\end{theo}

\begin{pf} We first notice that if $(\mathcal C, O)$ is in the Monge hypersuface $\mathcal M _\mathbb P$, then the projected arc from $O$ belongs to the Monge hypersurface $\mathcal M _{\mathbb P _O}$ where $\mathbb P _O$ is the $(n-1)$-dimensional projective space of lines through the point $O$. Indeed, this follows directly from Remark \ref{n+4}: recall the matrix $M$ from equation (\ref{arcmatrix}), and its submatrix $M_{00}$: the remark amounts to the existence of a basis in which the first $(n+4)$ columns of $M$ ($(n+4)$-truncation) take the form 
$$\left(\begin{array}{ccccccc}
1&0&\ldots &0&0&0&0\\
0&1&\ldots &0&0&0&0\\
0&0&\ldots &0&0&0&*\\
\vdots&\vdots&&\vdots&\vdots&\vdots&\vdots\\
0&0&\ldots & 0&0&0&*\\
0&0&\ldots & 1&0&0&0\\
\end{array}\right);$$
now it is clear that the $(n+3)$-truncation of $M_{00}$ leads, after an elementary transformation (between the second and last columns: this amounts to taking the series represented by the second line as a new parameter) to a matrix of the same shape. 

Using Theorem \ref{iota}, we deduce that the equation of $\mathcal M_\mathbb P$ is the image by $\iota _{n-1}$ of the equation of $\mathcal M _{\mathbb{P}_O}$.

In other words, the expression of this equation in terms of the coordinates $U_{\lambda}$ is independent of the dimension $n$, and we just have to compute it for $n=3$. This is a bit tricky, hence we produce a few details.

We already mentioned this case in Example \ref{conics}: the direct calculation is difficult to carry without a computer, thus we follow Halphen (see \cite{H3, H4, H5}). Due to the fact that the Monge invariant is of reduced bidegree $(3,3)$, we understand that it can be written as a linear combination of the following monomials: $$U^2U_{111}, U^2U_{21}, U^2U_3, UU_2, UU_{11}, U_1^3$$. 

In the chart $\mathcal U_{01}$ of $\mathcal A _{\mathbb{P}_{2}}$, the arc is defined by the triple of series $(1,t, \sum_{i\geq 2}a_it^i)$. The non-vanishing coordinates $U_\lambda$ (with $\vert \lambda \vert \leq 3$) are reduced to $U_i=a_{i+2}$ and $U=a_2$. The condition that the corresponding truncated arc of order $5$ is contained in a conic is the vanishing of the determinant
$$\mu = \det \left (\begin{array}{ccc}0&a_2&2a_3\\a_2&a_3&a_4\\a_3&a_4&a_5\\ \end{array}\right )=- 2a_3 ^3+3a_2a_3a_4 -a_2^2a_5$$
which is, up to sign, the expression of the theorem.

To find the terms involving a partition of length $2$, we go to the case $n=3$ and begin by the remark that the vector space of differential invariants of order $6$ and reduced degree $(3,3)$ over $\mathbb{P}_{3}$ is $1$-dimensional, since clearly $\mathcal{A}_{\mathbb{P}_{3},6}$ has a dense orbit under $SL(A\spcheck)$. The Monge hypersurface may therefore, as announced in Example \ref{conics}, be defined as  the wronskian of the tangent scroll. We compute it at the arc $(1,t,t^{2}+\Sigma_{4\leq i}\alpha_{i}t^{i},t^{3}+\Sigma_{4\leq i}\beta_{i}t^{i})$, for which $U=1, U_{1}=\beta_{4}, U_{2}=\beta_{5}, U_{11}=-\alpha_{4}, U_{3}=\beta_{6}, U_{21}=-\alpha_{5}, U_{111}=0$. We find
$$\left\vert\begin{array}{lcccccccccr}3&4\beta_{4}&5\beta_{5}&6\beta_{6}\\1&0&3\alpha_{4}&4\alpha_{5}\\0&2&3\beta_{4}&4\beta_{5}\\0&0&1&2\beta_{4}\end{array}\right\vert=12(-2\beta_{4}^{3}+3\beta_{4}\beta_{5}-3\alpha_{4}\beta_{4}+2\alpha_{5}-\beta_{6})$$
{\noi}which is, up to scalar, the expression of the theorem.

To find the missing coefficient of $U_{111}$ we follow a trick, of independent interest, explained in \cite{H5}, chap. III, \S 5. Let $a=\Sigma_{0\leq i}a_{i}t^{i}$ be a smooth parametrised arc in $\mathbb{P}_{2}$, consider its associated differential equation $Ux^{(3)}-3U_{1}x^{(2)}+6U_{11}x'-6U_{111}x=0$ as in Proposition \ref{equadif}. In order to compute the (invariant) Monge equation one may (by replacing the unknown power series $x$ by a suitable multiple by a given power series) assume that the coefficient $U_{1}$ of $x^{2}$ is $0$ (Liouville), and divide by $U$: we get the differential equation $x^{3}+6u_{11}x'-6u_{111}x=0$, whose adjoint is $y^{3}+6v_{11}y'-6v_{111}y=0$, with $v_{11}=u_{11}$, $v_{111}=-u'_{11}-u_{111}$.  

We already mentioned, in Proposition \ref{duality1}, that the arcs associated to a differential equation and to its adjoint are a pair of Pl\"{u}cker dual curves: this remark provides a {\it (Pl\"ucker) duality} involution of the ring $\mathcal{H}_{n}/(U-1)$ (the differential invariants, specialised at $U=1$: we cannot escape this specialisation once we work in the setting of differential equations), preserving the remaining graduation. Therefore it appears interesting to compute this involution.

We content ourselves by the straightforward calculation of the matrix of this involution in the stable vector space $\langle u_{3},u_{21},u_{111}\rangle$ which is 
$$\left(\begin{array}{lcccccr}-2&3/4&-3/10\\-4&3&-2/5\\0&0&1\end{array}\right).$$
 Since the Monge invariant is (assuming $u_{1}=0$) an eigenvector of this involution, we find $u_{3}-2u_{21}-2u_{111}$ for its expression, as announced. 

\end{pf}

\end{document}